\documentclass[12pt]{article}
\usepackage[utf8]{inputenc}
\usepackage{amsmath,amssymb,lscape,longtable,enumerate}
\usepackage[english,russian]{babel}
\usepackage{indentfirst}
\usepackage[dvips]{graphicx}
\usepackage{graphicx}
\usepackage{amsfonts}
\usepackage{amsthm}
\usepackage{mathtools}
\usepackage{systeme}
\usepackage{physics}
\usepackage[colorlinks=true, allcolors=blue]{hyperref}

\renewcommand{\le}{\leqslant}
\renewcommand{\ge}{\geqslant}
\renewcommand{\leq}{\leqslant}
\renewcommand{\geq}{\geqslant}
\begin{document}

\title{On upper bounds on the number of parts in the problem of partitioning sets into parts of smaller diameter}

\author{A. I. Bikeev\footnote{Moscow Institute of Physics and Technology (State University), Department of Discrete Mathematics}, A. M. Raigorodskii\footnote{Moscow Institute of Physics and Technology (State University), Department of Discrete Mathematics and Laboratory of Advanced Combinatorics and Network Applications;
Lomonosov Moscow State University, Faculty of Mechanics and Mathematics, Department of Mathematical Statistics and Random Processes; Caucasian Mathematical Center of Adyghe State University;
Buryat State University, Institute of Mathematics and Informatics.}}
\date{}
\maketitle

\section{Introduction}

In the present paper, we study problems related to the classical Borsuk's problem. Recall that the Borsuk's problem (see \cite{Bor}) consists of finding the smallest number $ f(n) $ of parts of smaller diameter into which an arbitrary set of diameter 1 in Euclidean space $ {\mathbb R}^n $ can be divided. Here we will discuss the quantity $ \chi(n,b) $, which differs from the quantity $ f(n) $ in that, in its definition, an arbitrary set of diameter 1 in $ {\mathbb R}^n $ must be partitioned not simply into parts of smaller diameter, but into parts whose diameters are strictly less than a given number $ b \in (0,1] $. Thus, $ \chi(n,1) = f(n) $. 

The function $ \chi(n,b) $ and some related functions have been studied in numerous papers. We will cite some of these papers below. In this note, we concentrate on the case where $ b $ is a fixed number (there are cases of non-constant $ b = b(n) $), аnd $ n \to \infty $. In this case, for $ b = 1 $, the best estimates are
$$
\left(\left(\frac{2}{\sqrt{3}}\right)^{\sqrt{2}}+o(1)\right)^{\sqrt{n}} = (1.2255\ldots + o(1))^{\sqrt{n}} \le \chi(n,1) = f(n) \le \left(\sqrt{\frac{3}{2}}+o(1)\right)^n = (1.224\ldots+o(1))^n.
\eqno{(1)}
$$
The lower bound was proved in \cite{Rai1}, and the upper bound was obtained in \cite{Sch} and \cite{BL} in two different ways. Let us immediately note that the notation of $ o(1) $ in the inequalities means that there exists a function of the argument $ n $, such that this function tends to zero as $ n \to \infty $, and the corresponding inequality remains true for any $n \in \mathbb{N} $ when replacing $o(1)$ with this function. Naturally, the corresponding function can
also take negative values in the lower bound and also positive values in the upper bound. There is a fairly extensive literature discussing explicit expressions for the quantities $f(n)$ and $\chi(n,b)$ (see, for example, the surveys \cite{Rai2}--\cite{Rai6} and the book \cite{Rai7}). However, this is relevant only for particular ``small'' values of $ n $ (as follows from inequalities (1), we do not even know the asymptotics of the logarithm of the function $f(n)$), and, as we already said, we are now concerned only with the behavior of $f(n)$ as $ n \to \infty $. 

If $ b < 1 $, then it is well known that the lower and upper bounds of $ \chi(n,b) $ are exponential with respect to $ n $. This shows how unstable the formulation of the classical Borsuk's problem is in a certain sense: while in the case $b=1$ not only the asymptotics of the logarithm of $\chi(n,b)$, but also the order of its growth (due to the gap between $\sqrt{n}$ in the lower bound and $n$ in the upper bound) is unknown, in the case $b=1$ only the asymptotics of the logarithm are unknown, and the order of its growth is $n$.

The optimal lower bounds are explicitly written in the paper \cite{Rai8} (see also \cite{Kon}). However, the upper bounds are scattered across different publications, and we do not know of a single source of information on them. In this paper, we collect such information together and, among other things, find a new upper bound for the quantity $ \chi(n,b) $. 

So, in the next section we carefully discuss the upper bounds for $ \chi(n,b) $ and formulate a new result, which we prove in the third section of the paper.

Before proceeding to the implementation of the stated plan, we note that, as we have already written above, there are a number of quantities closely related to the quantity $ \chi(n,b) $. For example, in the work \cite{RaiMin}, the following quantity was considered: the smallest number of parts into which an arbitrary set of diameter 1 in $ {\mathbb R}^n $ can be divided so that in each part there are no pairs of points at the distance $ b $ from each other. In other words, in this case we allow in each part not only distances smaller than $b$ but also greater than $b$. Also, a very significant literature is devoted to quantities that are in some sense inverse functions of $ \chi(n,b) $, namely, the quantities
$$
d_k^n  = \sup_{F \subset {\mathbb R}^n: \, {\rm diam}\, F = 1} \inf_{F_1, \dots, F_k: \, F = F_1 \cup \ldots \cup F_k} \max_{i = 1, \dots, k} {\rm diam}\,F_i.
$$
Thus, any upper bound of $ d_k^n $ means that for any set of diameter 1 in $ {\mathbb R}^n $ there exists a partitioning $ F $ into $ k $ parts of diameter less than or equal to this upper bound. In other words, in the definition of $ \chi(n,b) $ we minimize the number of parts of a fixed diameter, but in the definition of $ d_k^n $ we minimize the diameter of each part in the partitioning into a fixed number of parts. Obviously, the estimates of the quantities $ \chi(n,b) $ and $ d_k^n $ can be easily turned to each other. However, almost all studies of the quantities $ d_k^n $ have focused on cases of fixed dimensions (see \cite{Fil1}--\cite{Kit}), and it is impossible to deduce the consequences we need from their results.
Only one work \cite{Fil2} stands out in this series and we will discuss its result in the next section.

\section{Various upper bounds of $ \chi(n,b) $ for $ b \in (0,1) $}

Firstly, according the classical Jung's theorem (see \cite{Jung}, \cite{Hel}), every set of diameter 1 in $ {\mathbb R}^n $ can be covered by a ball (let us call it {\it Jung's ball}) of radius
$$
\sqrt{\frac{n}{2n+2}} \sim \frac{1}{\sqrt{2}}.
$$
By dividing this ball by pairwise orthogonal hyperplanes passing through the center, one can get $ 2^n $ parts of diameter strictly smaller than 1. This immediately gives the estimate
$$
\chi(n,1) = f(n) \le 2^n.
\eqno{(2)}
$$
The asymptotics of the estimate (2) are worse than the asymptotics of the estimate (1), but (2) is better than (1) in small dimensions when $o(1)$ takes any concrete values. Since we do not talk about small dimensions in this note, we will not further develop the idea of the boundary (2), although this is possible (see \cite{Kit}, \cite{Las}).

For increasing dimension and arbitrary $ b \in (0,1) $, the idea of covering the Jung's ball with balls of diameter $ b $ is much stronger than the idea with hyperplanes. In the classical article by Rogers (see \cite{Rog}) it was proved that the sphere of radius $ r $ can be covered by spherical ``caps'' of radius $ \rho < r $, using no more than

$$
\left(\frac{r}{\rho} + o(1)\right)^n 
\eqno{(3)}
$$
caps. A huge number of works are related to refinement of the infinitesimal term in the base of the exponent (see, for example, \cite{Ver}, \cite{BW}). But we will not write out their results here either, since we are interested exclusively in the asymptotics of the logarithm of the estimate. In our problem, it is necessary to cover not a sphere, but a ball, and for small $ b $ it is not at all obvious that this can be done. However, in the article \cite{Ver}, the estimate (3) is also obtained in the case of a ball (the authors used the natural idea of dividing the ball into spherical layers). Thus, we can cover the Jung's ball of radius $ \frac{1}{\sqrt{2}} + o(1) $ by at most
$$
\left(\frac{\sqrt{2}}{b}+o(1)\right)^n.
$$
balls of radius $ \frac{b}{2} $.
And then we can easily divide each of the covering balls into $ n+1 $ parts of strictly smaller diameter (see, for example, \cite{Rai3}), and this gives the first estimate of the type we are interested in:
$$
\chi(n,b) \le \left(\frac{\sqrt{2}}{b}+o(1)\right)^n.
\eqno{(4)}
$$

We have not encountered the estimate (4) in this form in the literature, strangely enough. This estimate can be improved using the result of Bourgain--Lindenstrauss from the paper \cite{BL}. The upper bound in the chain of inequalities (1) is proved in that paper. Namely, the authors added a number of non-trivial arguments to the idea of applying the Rogers bound (3) to the Jung's ball, which allows to prove that, in fact, each set of diameter 1 in $ {\mathbb R}^n $ can be covered by
$$
\left(\sqrt{\frac{3}{2}}+o(1)\right)^n = (1.224\ldots+o(1))^n
$$
balls of the same diameter. The technique we used for the bound (4) can be applied to each of the mentioned covering balls to get 
$$
\chi(n,b) \le \left(\sqrt{\frac{3}{2}}+o(1)\right)^n \cdot \left(\frac{1}{b}+o(1)\right)^n \cdot (n+1) = 
\left(\sqrt{\frac{3}{2}} \cdot \frac{1}{b} +o(1)\right)^n.
\eqno{(5)}
$$
For any $ b $, the inequality (5) has stronger asymptotics than the inequality (4). 

So far, we have covered sets by balls.
Another idea consists in the covering of a set by its own congruent copies or with the congruent copies obtained from it by reflection with respect to the origin and compression. Using this approach and the estimates (3) and (4) from the paper  \cite{RZ} of Rogers and Zong, we get the inequality
$$
\chi(n,b) \le \left(\frac{1}{b}+1+o(1)\right)^n.
\eqno{(6)}
$$  
Clearly, when the parameter $ b $ passes the boundary $ \sqrt{\frac{3}{2}} - 1 = 0.224\ldots $, the estimates (5) and (6) replace each other in terms of optimality. 

As we noted at the end of the introduction, the quantity $ \chi(n,b) $ in the mode we are interested in was discussed in the paper \cite{Fil2} (see section 4 of that paper). Unfortunately, the author does not refer to estimates (4)--(6) in any way. His result has the form
$$
\chi(n,b) \le \left(\sqrt{\frac{\pi e}{6}} \cdot \left(2+\frac{1}{b}\right)+o(1)\right)^n.
\eqno{(7)}
$$  
It is sometimes weaker and sometimes stronger than the estimate (5). However, it is always worse than the estimate (6).

We have not encountered any other explicitly written estimates of the type (4)--(7) in the literature. 
However, the estimates (6) and (7) are based on the same general idea that one can cover (possibly, with overlaps) the whole space by some (equivalent) sets $ B_i $ of diameter $ b $, and then show using the density of the covering that any given set can be placed so that it intersect sufficiently small number of $ B_i $.
This idea was used by Rogers and Zong (with congruent copies of a given body) and by Filimonov (with the permutohedrons). Clearly, one can use this approach and find more convenient $ B_i $. For example, we can apply the results related to the packings and coverings of the balls in the space. Some ideas from \cite{LR} and \cite{Pr} can be easily applied to obtain the inequality 
$$
\chi(n,b) \le \left(\frac{2}{b}+4+o(1)\right)^n,
\eqno{(8)}
$$  
which is obviously weaker than all the inequalities (4)-(7). So, in this approach, we were unable to achieve a stronger result than the record-breaking scores (5) and (6) mentioned above. However, we refined the estimate (5).

\vskip+0.2cm
  
\noindent {\bf Theorem 1.} {\it For any $ b \in (0,1] $, the following holds:
$$
\chi(n,b) \le \left(\sqrt{\frac{1}{b^2}+\frac{1}{2}}+o(1)\right)^n.
\eqno{(9)}
$$} 

\vskip+0.2cm

In fact, theorem 1 is a corollary of reworking the technique from the paper \cite{BL} for the case where we cover a set of diameter 1 in $ \mathbb{R}^n $ by balls of diameter $ b $ instead of balls of diameter 1. This gives an advantage over the bound (5), which directly uses the result from \cite{BL} and then performs a subpartition into sets of diameter $ b $. In principle, we do not exclude that something similar has been done elsewhere. However, we have not found such references and, for completeness, we present the proof of Theorem 1 in the next section.

The bound (9) is stronger than both bounds (5) and (6) which (as we saw above) are the best among estimates (4)--(7) of the type we needed. Indeed, the estimate (5), as expected, coincides with the estimate (9) for $b = 1$, and, obviously, the inequality
$$
\sqrt{\frac{1}{b^2}+\frac{1}{2}} < \sqrt{\frac{3}{2}} \cdot \frac{1}{b}
$$
holds for any $ b < 1 $. At the same time,
$$
\sqrt{\frac{1}{b^2}+\frac{1}{2}} < \frac{1}{b}+1 \Longleftrightarrow 
\frac{1}{b^2}+\frac{1}{2} < \frac{1}{b^2} + \frac{2}{b}+1 \Longleftrightarrow 
\frac{2}{b} > - \frac{1}{2},
$$
and the last condition is trivial.
Note that even in the asymptotics with respect to $ b \to 0 $ the value of the estimate from (5) is significantly weaker than the values of estimates from (6) and (9), because it is asymptotically greater than them $ \sqrt{\frac{3}{2}} = 1.224\ldots $ times. However, the bounds (6) and (9) asymptotically behave as $ \frac{1}{b} + 1 $ and $ \frac{1}{b} + (1+o(1))\frac{b}{4} $, respectively, that is, they differ only in the second term.

\section{Proof of Theorem 1}

The proof we present in the following is largely repeated from the argument in the paper \cite{BL}. Denote by $ D $ the number $ \frac{1}{b} $. We need to prove for any $ D \ge 1 $ the estimate
$$
\chi(n,b) = \chi(n,1/D) \le \left(\sqrt{D^2 + \frac{1}{2}} + o(1)\right)^n.
$$

The proof of the following lemma is actually in the paper \cite{BL} (there, as in ours, it is lemma 1), so we do not present it. The formulation in the paper \cite{BL} differs in the requirement $ \frac{5}{9} \leq r \leq 1$ instead of the our requirement $ \lambda D \leq r \leq D $  (see below). However, the requirements are consistent for $ D = 1 $, $ \lambda = \frac{5}{9} $, and practically the same proof remains true in this formulation. 
To make it easier to understand the meaning of the formulation of lemma 1 below, we immediately note that we need to cover sets of diameter $1$ with sets of diameter $ b $, namely, with spherical caps. And the desired covering is equivalent to the covering of a set of diameter $ D $ with sets of diameter 1. 

\vskip+0.2cm

\noindent {\bf Lemma 1.} {\it Let 
$$
D \geq 1, ~~ \lambda \in \left(\frac{1}{2}, \frac{5}{9}\right], ~~  \lambda D \leq r \leq D, ~~ \varepsilon > 0.
$$
Then for some $\delta = \delta(\varepsilon)>0, $ $k_0 = k_0(\varepsilon) \in {\mathbb N}$ and for any sufficiently large $n > n(\varepsilon)$ there exists a family of coverings
$$
\mathcal{A}_k = \left\{\{C_{k,m}\}_{m=1}^{M_k}, ~ k=0,1,2,\dots,k_0\right\}
$$
of the sphere $S_r^{n-1}$ with spherical caps $C_{k,m}$ of diameters $2\rho_k = (1+\delta)^k$ such that
$$
M_k \leq ((r+\varepsilon)/\rho_k)^n,\; 0 \leq k \leq k_0,
$$
$$
(1+\delta)^{k_0} \geq \frac{r}{\lambda},
$$
and for any $k \in \{1, \ldots, k_0\}$, $m \in \{1, \ldots, M_k\}$ the following holds:
$$
C_{k,m} \subset \bigcup_{j \in B_{k,m}} C_{k-1, j},\; |B_{k,m}| \leq (1+\varepsilon)^n.
$$}

\vskip+0.2cm

Note that in \cite{BL} the proof of the analogue of lemma 1 is given without detailed calculations. In principle, it would be possible to fill this gap here by proving lemma 1 with all the details. However, these details are of an extremely technical nature and we will not clutter the current note with them.

The following lemma generalizes a number of statements from \cite{BL}, but we present its proof in full, since it cannot be deduced from the work \cite{BL} immediately. However, this is a technical lemma. Also, for clarity, we note that the bound $ \frac{D}{\sqrt{2}} $ in the formulation is a value slightly larger than the radius of the Jung's ball, which we discussed at the beginning of Section 2. In the proof of Lemma 2, this bound will not be needed. We are simply not interested in large radii, since we know that every set of diameter 1 can be covered by the Jung's ball.

\vskip+0.2cm

\noindent {\bf Lemma 2.} {\it Denote
$$
\widetilde{\alpha}_D =\sqrt{\frac{D^2}{4}-\frac{1}{8}}.
$$
Let
$$
D \geq 1, ~~ 0 \le a \le \widetilde{\alpha}_D, ~~ \sqrt{a^2 + \frac{1}{4}} \leq r \leq \frac{D}{\sqrt 2}, ~~ \frac{1}{2} \leq \rho \leq r,
\eqno{(10)}
$$
and let
$$
f(r,\rho,a) =\frac{r}{\rho}\cdot\frac{a\sqrt{r^2 - \rho^2} + \sqrt{r^2 - \frac{1}{4}-a^2} \cdot \sqrt{\rho^2 - \frac{1}{4}}}{r^2 - \frac{1}{4}}.
$$
Then $f(r,\rho,a) \leq \sqrt{4a^2+1}$.}

\vskip+0.2cm

Note that 
$$
\frac{D}{\sqrt{2}} = \sqrt{2 \left(\widetilde{\alpha}_D\right)^2 + \frac{1}{4}},
$$
thus the segment of possible values of $r$ in the condition (10) is not empty even for
 $ a = \widetilde{\alpha}_D $. Moreover, it is not empty for smaller values of $a$.

\paragraph{Proof of Lemma 2.}  Consider the derivative of the function $ f $ with respect to $ \rho $:  
$$
\frac{\partial f(r,\rho, a)}{\partial\rho} = \frac{r}{r^2-\frac{1}{4}} \cdot \frac{1}{\rho^2 \sqrt{r^2-\rho^2}\sqrt{\rho^2-\frac{1}{4}}} \cdot \left( \frac{1}{4}\sqrt{r^2-\frac{1}{4}-a^2}\cdot\sqrt{r^2-\rho^2}-ar^2\sqrt{\rho^2-\frac{1}{4}} \right).$$

From the form of the derivative it is clear that, for fixed values of $r$ and $a$, the maximum of the function $f(r,\rho,a)$ on the interval $\rho\in \left[\frac{1}{2},r\right]$ is achieved under the condition
$$
\rho^2 = \frac{4a^2r^4 + \left(r^2-\frac{1}{4}-a^2\right)r^2}{16a^2r^4+r^2-\frac{1}{4}-a^2} = \frac{4a^2+1}{16a^2r^2+4a^2+1}\cdot r^2;
$$
$$
\frac{1}{\rho^2} = \frac{1}{r^2}+\frac{16a^2}{4a^2+1}.
$$

We are only left to calculate the value of this maximum. Let
$$
\frac{1}{\rho^2} = \frac{1}{r^2}+\frac{16a^2}{4a^2+1}, ~~ r^2 = a^2 + \frac{1}{4} + s, ~~ s \geq 0. 
$$
Then
$$
f(r,\rho,a) =\frac{r}{\rho}\cdot\frac{a\sqrt{r^2 - \rho^2} + \sqrt{r^2 - \frac{1}{4}-a^2} \cdot \sqrt{\rho^2 - \frac{1}{4}}}{r^2 - \frac{1}{4}} = 
$$
$$
=r \cdot \frac{ar\sqrt{\frac{1}{\rho^2}-\frac{1}{r^2}}+\frac{\sqrt{s}}{2}\sqrt{4-\frac{1}{\rho^2}}}{a^2+s} = \frac{r}{a^2+s} \cdot \left(ar\sqrt{\frac{16a^2}{4a^2+1}}+\frac{\sqrt{s}}{2}\sqrt{\frac{4}{4a^2+1}-\frac{1}{r^2}}\right)=
$$
$$
=\frac{r}{(a^2+s)\sqrt{4a^2+1}}\cdot\left(4a^2r+\frac{\sqrt{s}}{r}\sqrt{r^2-a^2-\frac{1}{4}}\right) = \frac{r}{(a^2+s)\sqrt{4a^2+1}}\cdot \left(4a^2 r+\frac{s}{r}\right) =
$$
$$
=\frac{4a^4+4a^2s+a^2+s}{(a^2+s)\sqrt{4a^2+1}} = \frac{4a^2+1}{\sqrt{4a^2+1}} = \sqrt{4a^2+1}.
$$

The proof of lemma 2 is complete.

\vskip+0.2cm

The following lemma is a generalization of lemma 2 from the paper \cite{BL}. However, the generalization is not obvious, and we provide a detailed proof below.

\vskip+0.2cm

\noindent {\bf Lemma 3.} {\it Let
$$
D \ge 1, ~~ \sqrt{\left(\widetilde{\alpha}_D\right)^2 + \frac{1}{4}} \le r \le \sqrt{2 \left(\widetilde{\alpha}_D\right)^2 + \frac{1}{4}} = \frac{D}{\sqrt{2}}, ~~ \frac{1}{2} \leq \rho \leq r.
$$
Let $K$ be a subset of diameter $D$ on the sphere $S_r^{n-1}$. Let $C_1$ and $C_2$ be two spherical caps of diameter $2\rho$ on the sphere $S_r^{n-1}$. Suppose that each of the sets $C_1 \cap K$ и $C_2 \cap K$ is not contained in a union of two spherical caps of diameter $1$. Let the distance between the centers of $C_1$ and $C_2$ be $2d$. Then
$$
\frac{d}{\rho} \leq \sqrt{D^2 + \frac{1}{2}}.
$$}

\vskip+0.2cm

Note that the lower bound of $ r $ in lemma 3 is
$$
\sqrt{\frac{D^2}{4} + \frac{1}{8}} > \frac{D}{2}.
$$
This is useful to understand when talking about sets of diameter $ D $ on the sphere of radius $ r $.

\paragraph{Proof of lemma 3.} If the caps $C_1$ and $C_2$ have non-empty intersection, then the distance between the centers of $C_1$ and $C_2$ is at most $2\rho$. Then
$$
\frac{d}{\rho} \leq \frac{\rho}{\rho} = 1 < \sqrt{\frac{3}{2}} \leq \sqrt{D^2+\frac{1}{2}}.
$$
Therefore, we further assume that the intersection of $C_1$ and $C_2$ is not empty.

Let us consider a coordinate system such that the center of the cap $C_1$ has coordinates
$$
\left(\sqrt{r^2 - d^2}, d, 0,0,\dots\right), 
$$ 
and the center of the cap $C_2$ has coordinates
$$
\left(\sqrt{r^2 - d^2}, -d, 0,0,\dots\right). 
$$
For any
$$
\alpha \in \left[0,\sqrt{r^2-\frac{1}{4}}\right] ,
$$ 
let us consider the ball $B_{\alpha} = B(q_1, 1/2)$ with center 
$$
q_1 = \left(\sqrt{r^2 - \frac{1}{4} - \alpha^2}, \alpha, 0,0,\dots\right)
$$ 
and radius $\frac{1}{2}$. Let $\alpha$ be the maximal number such that any point of the set $C_1 \setminus B_{\alpha}$ has the second coordinate at least $\alpha$ (that is, $B_\alpha$ covers a ``bottom part'' of the cap $C_1$), see fig. 1.

\begin{figure}[ht]
\center{\includegraphics[scale=0.5, width=230pt]{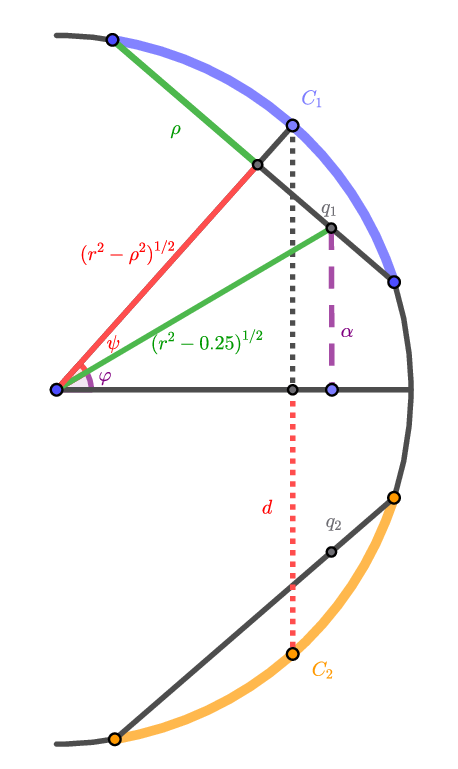}}
\end{figure}

Since $C_1 \cap K$ is not contained in a union of two balls of diameter $\frac{1}{2}$, the set $K \cap C_1 \setminus B_\alpha$ is not contained in a ball of diameter $1$. Thus, the closure of $K \cap C_1 \setminus B_\alpha$ contains a set $\{{\bf x}_j\}_{j=1}^m$ such that the ball $\Omega_1$ of minimal radius containing these points has a radius $t_1 > 1/2$ and such that all the points ${\bf x}_j$ are on the boundary of the ball $\Omega_1$. Then center ${\bf u}$ of the ball $\Omega_1$ can be represented as
$$
{\bf u} = \sum\limits_{j=1}^m \lambda_j {\bf x}_j, ~~  \lambda_j\geq 0, ~~ \sum\limits_{j=1}^m \lambda_j = 1,
$$
since ${\bf u}$ is contained in the convex hull of the points ${\bf x}_j$. The second coordinate of each ${\bf x}_j$ and therefore also of ${\bf u}$ is at least $\alpha$. We also have $|{\bf u}-{\bf x}_j|=t_1$ for any $j \in \{1, \dots, m\}$.

Similarly, we can find points $ {\bf v} $, $ \{{\bf y}_i\}_{i=1}^s $ in the closure of $K \cap C_2$ so that
$$
{\bf v}=\sum\limits_{i=1}^s \mu_i {\bf y}_i, ~~ \mu_i\geq 0,  ~~ \sum\limits_{i=1}^s \mu_i = 1,
$$
the second coordinate of is at most $-\alpha$ and so that $|{\bf v}-{\bf y}_i|=t_2 > \frac{1}{2}$ for any $i\in \{1, \dots, s\}$.

Thus, $|{\bf u}-{\bf v}| \geq 2\alpha$. Since
$$
|{\bf x}_j-{\bf y}_i|^2 = |{\bf x}_j - {\bf u}|^2 + |{\bf u}-{\bf v}|^2 + |{\bf v}-{\bf y}_i|^2 +
$$
$$
+ 2(<{\bf x}_j-{\bf u}, {\bf u}-{\bf v}> + <{\bf u}-{\bf v}, {\bf v}-{\bf y}_i> + <{\bf x}_j-{\bf u}, {\bf v}-{\bf y}_i>),
$$
we get
$$
\sum\limits_{i=1}^s\sum\limits_{j=1}^m \mu_i \lambda_j |{\bf x}_j-{\bf y}_i|^2 = |{\bf u}-{\bf v}|^2+t_1^2+t_2^2 > (2\alpha)^2 + 2\cdot(1/2)^2 = 4\alpha^2 + \frac{1}{2}.
$$

Since for any $i$ and $j$ the distance $|{\bf x}_j-{\bf y}_i|$ is at most $D$, we get
$$
\sum\limits_{i=1}^s\sum\limits_{j=1}^m \mu_i \lambda_j |{\bf x}_j-{\bf y}_i|^2 \leq D^2 .
$$ 
Hence,
$$
D^2 > 4\alpha^2 + \frac{1}{2} \Longrightarrow \alpha <  \sqrt{\frac{D^2}{4} - \frac{1}{8}} = \widetilde{\alpha}_D.
$$

Denote by $\varphi$ and $\psi$ corners as on fig. 1. It is easy to see that
$$
\sin \varphi = \frac{\alpha}{\sqrt{r^2 - \frac{1}{4}}}, ~~ \cos \varphi = \frac{\sqrt{r^2-\frac{1}{4}-\alpha^2}}{\sqrt{r^2 - \frac{1}{4}}},
$$
$$
\sin \psi = \frac{\sqrt{\rho^2 - \frac{1}{4}}}{\sqrt{r^2 - \frac{1}{4}}}; ~~ \cos \psi = \frac{\sqrt{r^2-\rho^2}}{\sqrt{r^2 - \frac{1}{4}}}.
$$

The distance between the centers of the caps $C_1$ и $C_2$ is $2d$. Since $d = r\sin(\varphi+\psi)$, we get
$$
d \leq r\cdot\frac{\alpha\sqrt{r^2 - \rho^2} + \sqrt{r^2 - \frac{1}{4}-\alpha^2} \cdot \sqrt{\rho^2 - \frac{1}{4}}}{r^2 - \frac{1}{4}}.
$$
Then due to lemma 2 we have
$$
\frac{d}{\rho} \leq f(r,\rho,\alpha) \leq \sqrt{4\alpha^2+1} \leq \sqrt{4\cdot\left(\frac{D^2}{4} - \frac{1}{8}\right)+1} = \sqrt{D^2+\frac{1}{2}}.
$$

The proof of lemma 3 is complete. 

\paragraph{Completion of the proof of Theorem 1.} Here everything is very similar to the work \cite{BL}, so we will give only a very brief outline.

Firstly, note that we deal with spheres instead of balls, but one can be reduced to other by inscribing a sufficiently large number of concentric spheres into a ball (cf. the remark to formula (3) in Section 2).

Next, we know that any set of diameter $ D $ is contained in a ball of radius $ r \le \frac{D}{\sqrt{2}} $. There are two possible cases:
$$
r \le \sqrt{\left(\widetilde{\alpha}_D\right)^2+\frac{1}{4}}
\eqno{(11)}
$$
and
$$
\sqrt{\left(\widetilde{\alpha}_D\right)^2+\frac{1}{4}} \le r \le \sqrt{2\left(\widetilde{\alpha}_D\right)^2+\frac{1}{4}} = \frac{D}{\sqrt{2}}.
\eqno{(12)}
$$
In the first case, we simply use the estimates of Rogers (see \cite{Rog}, \cite{Ver}) and obtain the necessary bound:
$$
\frac{\sqrt{\left(\widetilde{\alpha}_D\right)^2+\frac{1}{4}}}{1/2} = \sqrt{D^2+\frac{1}{2}}.
$$

In the case of inequalities (12), let us introduce the parameter
$$
\lambda = \min\left(\frac{5}{9}, \frac{\sqrt{\frac{D^2}{4} + \frac{1}{8}}}{D}\right) \in \left(\frac{1}{2}, \frac{5}{9}\right].
$$
It satisfies the condition of lemma 1, and the parameter $r$ is also contained in the necessary segment:
$$
r \ge \sqrt{\left(\widetilde{\alpha}_D\right)^2+\frac{1}{4}} =  \sqrt{\frac{D^2}{4}+\frac{1}{8}} \ge \lambda D, ~~ r \le \frac{D}{\sqrt{2}} < D,
$$
that is, lemma 1 can be applied. Lemma 3 also holds in the conditions (12). Note that lemma 1 is an analogue of lemma 1 from \cite{BL}, and lemma 3 is an analogue of lemma 2 from the same paper. Accordingly, we turn to the section {\bf Proof of the upper estimate in (1)} on page 142 of the paper \cite{BL}. The text of this proof can now be reproduced almost word for word, starting from its second paragraph. Our Lemma 2, however, is not there, since it serves as a non-trivial generalization of the arguments located directly in the text of the article \cite{BL} on page 143 in the paragraph between formulas (6) and (7). Accordingly, the quantity $ d(r,\rho_k) $, which appears in the article \cite{BL} and turns (after dividing by $ \rho_k $) into $ \sqrt{\frac{3}{2}} $ only in the mentioned paragraph, is immediately estimated by us as $ \rho_k \sqrt{D^2+\frac{1}{2}} $. Further, the value $ \frac{10}{9} $, which appears in the second string after the formula (7) on the page 143 and is estimated by $ \sqrt{\frac{3}{2}} $ with a large margin in \cite{BL}, in our case turns to $ 2\lambda $. This value $2\lambda$ is less or equal to the necessary threshold $ \sqrt{D^2+\frac{1}{2}} $ (also by a large margin, because we neglect the denominator $ D $). 

There are no other changes. Of course, in the argument from \cite{BL}, we need to change the various coefficients of the terms of order $ \varepsilon $ in the bases of the exponents. But this does not affect the final conclusion, in which all these terms are collected in $ o(1) $, and we will not write all these technical details here.

The proof of Theorem 1 is complete.

\end{document}